\documentclass[10pt]{article}
\usepackage{graphicx}
\usepackage{amssymb}
\newtheorem{theorem}{Theorem}

\newtheorem{lemma}[theorem]{Lemma}

\input tcilatex

\begin{document}

\begin{titlepage}

\vskip0.2truecm

\begin{center}

{\LARGE \bf A consequence of the growth of rotation sets for families of diffeomorphisms of the torus}

\end{center}

\vskip  0.4truecm

\centerline {{\large Salvador Addas-Zanata}}

\vskip 0.2truecm

\centerline { {\sl Instituto de Matem\'atica e Estat\'\i stica }}
\centerline {{\sl Universidade de S\~ao Paulo}}
\centerline {{\sl Rua do Mat\~ao 1010, Cidade Universit\'aria,}} 
\centerline {{\sl 05508-090 S\~ao Paulo, SP, Brazil}}
 
\vskip 0.7truecm

\begin{abstract}

In this paper we consider $C^\infty $-generic families of area-preserving 
diffeomorphisms of the torus homotopic to the identity and their rotation sets. 
Let $f_t:\rm{T^2\rightarrow T^2}$ be such a family, 
$\widetilde{f}_t:\rm I\negthinspace R^2 \rightarrow \rm I\negthinspace R^2$ 
be a fixed family of lifts and $\rho (\widetilde{f}_t)$ be their rotation sets, 
which we assume to have interior for $t$ in a certain open interval $I.$ 
We also assume that some rational point 
$(\frac pq,\frac lq)\in \partial \rho (\widetilde{f}_{\overline{t}})$ 
for a certain parameter $\overline{t}\in I$ and we want to understand 
consequences of the following hypothesis: For all $t>\overline{t},$ 
$t\in I,$ $(\frac pq,\frac lq)\in int(\rho (\widetilde{f}_t)).$

Under these very natural assumptions, we prove that there exists a 
$f_{\overline{t}}^q$-fixed hyperbolic 
saddle $P_{\overline{t}}$ such that its rotation vector 
is $(\frac pq,\frac lq)$ and,
there exists a sequence 
$t_i>\overline{t},$ $t_i\rightarrow \overline{t},$ such that if $P_t$ is the 
continuation of $P_{\overline{t}}$ with the parameter, then 
$W^u(\widetilde{P}_{t_i})$ (the unstable manifold) has quadratic tangencies with 
$W^s(\widetilde{P}_{t_i})+(c,d)$ (the stable manifold translated by $(c,d)),$ 
where $\widetilde{P}_{t_i}$ is any 
lift of $P_{t_i}$  to the plane, in other words, $\widetilde{P}_{t_i}$ is a fixed point 
for $(\widetilde{f}_{t_i})^q-(p,l),$ and $(c,d)\neq (0,0)$ are certain integer vectors such 
that $W^u(\widetilde{P}_{\overline{t}})$ do not intersect 
$W^s(\widetilde{P}_{\overline{t}})+(c,d).$ And these tangencies become
transverse as $t$ increases.

As we also proved that for $t>\overline{t},$ $W^u(\widetilde{P}_t)$ 
has transverse intersections with $W^s(\widetilde{P}_t)+(a,b),$ for all integer 
vectors $(a,b),$ one may consider that the tangencies above are associated to
the birth of the 
heteroclinic intersections in the plane that did not exist
for $t\leq \overline{t}.$

\end{abstract} 

\vskip 0.3truecm

\noindent{\bf e-mail:} sazanata@ime.usp.br

\vskip 0.3truecm


\vskip 0.3truecm

\vfill
\hrule
\noindent{\footnotesize{The author is partially supported 
by CNPq, grant: 306348/2015-2}}

\end{titlepage}

\baselineskip=6.2mm

\section{Introduction}

\subsection{General explanations}

In this paper, in a certain sense, we continue the study initiated in \cite
{eupatrice}. There we looked at the following problem: Suppose $f:{\rm %
T^2\rightarrow T^2}$ is a homeomorphism homotopic to the identity and its
rotation set, which is supposed to have interior, has a point $\rho $ in its
boundary with both coordinates rational. The question studied was the
following: is it possible to find two different arbitrarily small $C^0$%
-perturbations of $f,$ denoted $f_1$ and $f_2$ in a way that $\rho $ does
not belong to the rotation set of $f_1$ and $\rho $ is contained in the
interior of the rotation set of $f_2?$ In other words we were asking if the
rational mode locking found by A. de Carvalho, P. Boyland and T. Hall \cite
{andre} in their particular family of homeomorphisms was, in a certain
sense, a general phenomenon or not. Our main theorems and examples showed
that the answer to this question depends on the set of hypotheses assumed.
For instance, regarding $C^2$-generic families, we proved that if $\rho \in
\partial \rho (\widetilde{f}_{\overline{t}})$ for some $\overline{t}$ and
for $t<\overline{t},$ close to $\overline{t},$ $\rho \notin \rho (\widetilde{%
f}_t),$ then for all sufficiently small $t-\overline{t}>0,$ $\rho \notin
int(\rho (\widetilde{f}_t)).$

Now we are interested in the dynamical consequences of a situation, which
may be obtained as a continuation of the previous one: Suppose $f_t:{\rm %
T^2\rightarrow T^2}$ is a one parameter family of diffeomorphisms of the
torus homotopic to the identity, for which the rotation set $\rho ( 
\widetilde{f}_t)$ at a certain parameter $t=\overline{t},$ has interior,
some rational vector $\rho \in \partial \rho (\widetilde{f}_{\overline{t}})$
and for all sufficiently small $t>\overline{t},$ $\rho \in int(\rho ( 
\widetilde{f}_t))).$ We want to understand what happens for the family $f_t,$
$t>\overline{t}.$ In other words, we are assuming that at $t=\overline{t}$
the rotation set is ready to locally grow in a neighborhood of $\rho .$

This is the usual situation for (generic) families:\ As the parameter
changes, the rotation set hits a rational vector, this vector stays for a
while in the boundary of the rotation set and finally, it is eaten by the
rotation set, that is, it becomes an interior point.

Here we consider $C^\infty $-generic area preserving families in the sense
of Meyer \cite{meyer1970} and also satisfying other generic conditions, and
the theorem proved goes in the following direction: If for instance, $f_t:%
{\rm T^2\rightarrow T^2}$ is such a family for which the rotation set at $t= 
\overline{t}$ has interior, $(0,0)\in \partial \rho (\widetilde{f}_{
\overline{t}})$ and for all sufficiently small $t>\overline{t},$ $(0,0)\in
int(\rho (\widetilde{f}_t)),$ then $\widetilde{f}_{\overline{t}}:{\rm I%
\negthinspace R^2\rightarrow }$ ${\rm I\negthinspace R^2}$ has a hyperbolic
fixed saddle $\widetilde{P}_{\overline{t}}$, for which $W^u(\widetilde{P}%
_{t_i})$ ($\widetilde{P}_t$ is the continuation of $\widetilde{P}_{\overline{%
t}}$ with the parameter $t$ and $t_i>\overline{t}$ is a certain sequence
converging to $\overline{t}$) has heteroclinic tangencies with certain
special integer translates of the stable manifold, $W^s(\widetilde{P}%
_{t_i})+(a,b),$ $(a,b)\in {\rm Z\negthinspace
\negthinspace Z^2,}$ which unfold (=become transversal) as $t$ increases.
The integer vectors $(a,b)$ mentioned above belong to a set $K_{{\rm Z%
\negthinspace
\negthinspace Z^2}}\subset {\rm Z\negthinspace
\negthinspace Z^2}$ and satisfy the following: for $t\leq \overline{t},$ $%
W^u(\widetilde{P}_t)$ can not have intersections with $W^s(\widetilde{P}%
_{t})+(a,b)$ whenever $(a,b)\in K_{{\rm Z\negthinspace
\negthinspace Z^2}}.$ Moreover, the sequence $t_i\rightarrow \overline{t}$
depends on the choice of $(a,b)\in K_{{\rm Z\negthinspace
\negthinspace Z^2}}.$

So increasing the parameter until the critical value at $t=\overline{t}$ is
reached (the moment when the rotation set is ready to locally grow), this
critical parameter is accumulated from the other side by parameters at which
there are heteroclinic tangencies in the plane (homoclinic in the torus) not
allowed to exist when $t\leq \overline{t}.$ In fact, as we will prove, for
any $t>\overline{t},$ 
\begin{equation}
\label{allint}W^u(\widetilde{P}_t)\text{ has transverse intersections with }%
W^s(\widetilde{P}_t)+(a,b),\text{ }\forall (a,b)\in {\rm Z\negthinspace
\negthinspace Z^2}. 
\end{equation}
In this way, the creation of the heteroclinic intersections for integers $%
(a,b)$ in (\ref{allint}) which did not exist for $t\leq \overline{t},$ also
produces tangencies.

In order to state things clearly and to precisely present our main result, a
few definitions are necessary.

\subsection{Basic notation and some definitions}

\begin{enumerate}
\item  Let ${\rm T^2}={\rm I}\negthinspace {\rm R^2}/{\rm Z\negthinspace
\negthinspace Z^2}$ be the flat torus and let $p:{\rm I}\negthinspace {\rm %
R^2}\longrightarrow {\rm T^2}$ be the associated covering map. Coordinates
are denoted as $\widetilde{z}\in {\rm I}\negthinspace {\rm R^2}$ and $z\in 
{\rm T^2.}$

\item  Let $Diff_0^r({\rm T^2})$ be the set of $C^r$ diffeomorphisms $%
(r=0,1,...,\infty )$ of the torus homotopic to the identity and let $%
Diff_0^r({\rm I\negthinspace R^2})$ be the set of lifts of elements from $%
Diff_0^r({\rm T^2})$ to the plane. Maps from $Diff_0^r({\rm T^2})$ are
denoted $f$ and their lifts to the plane are denoted $\widetilde{f}.$

\item  Let $p_{1,2}:{\rm I}\negthinspace {\rm R^2}\longrightarrow {\rm I}%
\negthinspace {\rm R}$ be the standard projections, respectively in the
horizontal and vertical components;

\item  Given $f\in Diff_0^0({\rm T^2})$ (a homeomorphism) and a lift $
\widetilde{f}\in Diff_0^0({\rm I}\negthinspace {\rm R^2}),$ the so called
rotation set of $\widetilde{f},$ $\rho (\widetilde{f}),$ can be defined
following Misiurewicz and Ziemian \cite{misiu} as: 
\begin{equation}
\label{rotsetident}\rho (\widetilde{f})=\bigcap_{{\ 
\begin{array}{c}
i\geq 1 \\ 
\end{array}
}}\overline{\bigcup_{{\ 
\begin{array}{c}
n\geq i \\ 
\end{array}
}}\left\{ \frac{\widetilde{f}^n(\widetilde{z})-\widetilde{z}}n:\widetilde{z}%
\in {\rm I}\negthinspace {\rm R^2}\right\} }
\end{equation}

This set is a compact convex subset of ${\rm I\negthinspace R^2}$ (see \cite
{misiu}), and it was proved in \cite{franksrat} and \cite{misiu} that all
points in its interior are realized by compact $f$-invariant subsets of $%
{\rm T^2,}$ which can be chosen as periodic orbits in the rational case. By
saying that some vector $\rho \in \rho (\widetilde{f})$ is realized by a
compact $f$-invariant set, we mean that there exists a compact $f$-invariant
subset $K\subset {\rm T^2}$ such that for all $z\in K$ and any $\widetilde{z}%
\in p^{-1}(z)$ 
\begin{equation}
\label{deffrotvect}\stackunder{n\rightarrow \infty }{\lim }\frac{\widetilde{f%
}^n(\widetilde{z})-\widetilde{z}}n=\rho .
\end{equation}
Moreover, the above limit, whenever it exists, is called the rotation vector
of the point $z,$ denoted $\rho (z)$$.$
\end{enumerate}

\subsection{Some background and the main theorem}

\subsubsection{Prime ends compactification of open disks}

If $U\subset {\rm I\negthinspace R^2}$ is an open topological disk whose
boundary is a Jordan curve and $\widetilde{f}:{\rm I\negthinspace %
R^2\rightarrow I\negthinspace R^2}$ is an orientation preserving
homeomorphism such that $\widetilde{f}(U)=U,$ it is easy to see that $
\widetilde{f}:\partial U\rightarrow \partial U$ is conjugate to a
homeomorphism of the circle, and so a real number $\rho (U)=rotation$ $%
number $ $of$ $\widetilde{f}\mid _{\partial U}$ can be associated to this
problem. Clearly, if $\rho (U)$ is rational, there exists a periodic point
in $\partial U$ and if it is not, then there are no such points. This is
known since Poincar\'e. The difficulties arise when we do not assume $%
\partial U$ to be a Jordan curve.

The prime ends compactification is a way to attach to $U$ a circle called
the circle of prime ends of $U,$ obtaining a space $U\sqcup S^1$ with a
topology that makes it homeomorphic to the closed unit disk. If, as above we
assume the existence of a planar orientation preserving homeomorphism $
\widetilde{f}$ such that $\widetilde{f}(U)=U,$ then $\widetilde{f}\mid _U$
extends to $U\sqcup S^1.$ The prime ends rotation number of $\widetilde{f}%
\mid _U,$ still denoted $\rho (U),$ is the usual rotation number of the
orientation preserving homeomorphism induced on $S^1$ by the extension of $
\widetilde{f}\mid _U.$ But things may be quite different in this setting. In
full generality, it is not true that when $\rho (U)$ is rational, there are
periodic points in $\partial U$ and for some examples, $\rho (U)$ is
irrational and $\partial U$ is not periodic point free. Nevertheless, in the
area-preserving case which is the case considered in this paper, many
interesting results were obtained. We refer the reader to \cite{mather}, 
\cite{frankslecal}, \cite{patmeykoro1} and \cite{patmeykoro2}. To conclude,
we present some results extracted from these works, adapted to our hypotheses.

Assume $h:{\rm T^2\rightarrow T^2}$ is an area-preserving diffeomorphism of
the torus homotopic to the identity such that for each integer $n>0,$ $h$
has finitely many $n$-periodic points. Moreover, we also assume more
technical conditions on $h:$ for each $n>0,$ at all $n$-periodic points a
Lojasiewicz condition is satisfied, see \cite{dumortier}. And if the
eigenvalues of $Dh^n$ at such a periodic point are both equal to $1,$ then
the point is topologically degenerate; it has zero topological index. As
explained in section 2 of \cite{eupatrice}, the dynamics near such a point
is similar to the one in figure 2. In particular $h$ has one stable
separatrix (like a branch of a hyperbolic saddle) and an unstable one at
such a periodic point, both $h$-invariant. Topologically, the local dynamics
in a neighborhood of the periodic point is obtained by gluing exactly two
hyperbolic sectors.

Fix some $\widetilde{h}:{\rm I\negthinspace %
R^2\rightarrow I\negthinspace R^2,}$ a lift of $h$ to the plane.
Given a $\widetilde{h}$-invariant continuum $K\subset {\rm I\negthinspace %
R^2,}$ if $O$ is a connected component of $K^c$ ($O$ is a topological open
disk in the sphere $S^2\stackrel{def.}{=}{\rm I\negthinspace R^2}\sqcup
\infty ,$ the one point compactification of the plane, that is, $O$ is a
connected simply connected open subset of $S^2$), which is also assumed to
be $\widetilde{h}$-invariant (it could be periodic with period larger than
1), let $\alpha $ be the rotation number of the prime ends compactification
of $O.$ From the hypothesis on $h,$ we have:
    
\vskip0.2truecm

{\bf Theorem A.} {\it If $\alpha $ is rational, then $\partial O$ has
accessible $\widetilde{h}$-periodic points. And if such a point has period $%
n,$ the eigenvalues of $D\widetilde{h}^n$ at this periodic point must be
real and can not be equal to $-1.$ So from the above properties assumed on $%
h,$ in $\partial O$ we either have accessible hyperbolic periodic saddles or
periodic points with both eigenvalues equal to $1,$ whose local dynamics is
as in figure 2. And then, there exist connections between separatrices of
the periodic points, these separatrices being either stable or unstable
branches of hyperbolic saddles, or the unstable or the stable separatrix of
a point as in figure 2.}

\vskip0.2truecm

The existence of accessible $\widetilde{h}$-periodic points can be found in 
\cite{cartlit}. The information about the eigenvalues is a new result from 
\cite{patmeykoro2} and the existence of connections in the above situation
can be found in \cite{mather}, \cite{frankslecal} and also \cite{patmeykoro2}.

\vskip0.2truecm

{\bf Theorem B.} {\it If $\alpha $ is irrational and $O$ is bounded, then
there is no periodic point in $\partial O.$}

\vskip0.2truecm

This is a result from \cite{patmeykoro1}.

\subsubsection{Some results on diffeomorphisms of the torus homotopic to the
identity}

As the rotation set of a homeomorphism of the torus homotopic to the
identity is a compact convex subset of the plane, there are three
possibilities for its shape:

\begin{enumerate}
\item  it is a point;

\item  it is a linear segment;

\item  it has interior;
\end{enumerate}

We consider the situation when the rotation set has interior.

Whenever a rational vector $(p/q,l/q)\in int(\rho (\widetilde{f}))$ for some 
$\widetilde{f}\in Diff_0^2({\rm I\negthinspace R^2}),$%
\begin{equation}
\label{resvida} 
\begin{array}{c}
\widetilde{f}^q(\bullet )-(p,l)\text{ has a hyperbolic periodic saddle } 
\widetilde{P}\in {\rm I\negthinspace R^2} \\ \text{such that }W^u(\widetilde{%
P}),\text{ its unstable manifold, has a topologically transverse} \\ \text{%
intersection with }W^s(\widetilde{P})+(a,b),\ \text{for all integer vectors }%
(a,b). 
\end{array}
\end{equation}
That is, the unstable manifold of $\widetilde{P}$ intersects all integer
translations of its stable manifold$.$ This result is proved in \cite
{c1epsilon}. Now it is time to precisely define what a topologically
transverse intersection is:

\begin{description}
\item[Definition (Top. Trans. Intersections):] If $f:M\rightarrow M$ is a $C^1$ diffeomorphism of an
orientable boudaryless surface $M$ and $p,q\in M$ are $f$-periodic saddle
points, then we say that $W^u(p)$ has a topologically transverse
intersection with $W^s(q),$ whenever there exists a point $z\in W^s(q)\cap
W^u(p)$ ($z$ clearly can be chosen arbitrarily close to $q$ or to $p$) and
an open topological disk $B$ centered at $z,$ such that $B\backslash \alpha
=B_1\cup B_2,$ where $\alpha $ is the connected component of $W^s(q)\cap B$
which contains $z,$ with the following property: there exists a closed
connected piece of $W^u(p)$ denoted $\beta $ such that $\beta \subset B,$ $%
z\in \beta ,$ and $\beta \backslash z$ has two connected components, one
contained in $B_1\cup \alpha $ and the other contained in $B_2\cup \alpha ,$
such that $\beta \cap B_1\neq \emptyset $ and $\beta \cap B_2\neq \emptyset .
$ Clearly a $C^1$ transverse intersection is topologically transverse. See
figure 1 for a sketch of some possibilities. Note that as $\beta \cap \alpha 
$ may contain a connected arc containing $z,$ the disk $B$ may not be chosen
arbitrarily small.

In order to have a picture in mind, consider $z$ close to $q,$ so that $z$
belongs to a connected arc in $W^s(q)$ containing $q,$ which is almost a
linear segment. Therefore, it is easy to find $B$ as stated above, it could
be chosen as an Euclidean open ball. Clearly, this is a symmetric definition: we can
consider a negative iterate of $z,$ for some $n<0$ such that $f^n(z)$
belongs to a connected piece of $W^u(p)$ containing $p,$ which is also
almost a linear segment. Then, a completely analogous construction can be
made, switching stable manifold with unstable: choose an open Euclidean ball 
$B^{\prime }$ centered at $f^n(z),$ such that $B^{\prime }\backslash \beta
^{\prime }=B_1^{\prime }\cup B_2^{\prime },$ where $\beta ^{\prime }$ is the
connected component of $W^u(p)\cap B^{\prime }$ which contains $z,$ with the
following property: there exists a closed connected piece of $W^s(q)$
denoted $\alpha ^{\prime }$ such that $\alpha ^{\prime }$$\subset B^{\prime
},$ $f^n(z)\in \alpha ^{\prime },$ and $\alpha ^{\prime }$$\backslash f^n(z)$
has two connected components, one contained in $B_1^{\prime }\cup \beta
^{\prime }$ and the other contained in $B_2^{\prime }\cup \beta ^{\prime },$
such that $\alpha ^{\prime }\cap B_1^{\prime }\neq \emptyset $ and
$\alpha ^{\prime }\cap B_2^{\prime }\neq \emptyset .$ So, $f^{-n}(B^{\prime }),$ $%
f^{-n}(\beta ^{\prime })$ and $f^{-n}(\alpha ^{\prime })$ are the
corresponding sets at $z.$
\end{description}

\vskip 0.2truecm

\begin{figure}[ht!]

\hfill

\includegraphics [height=78mm]{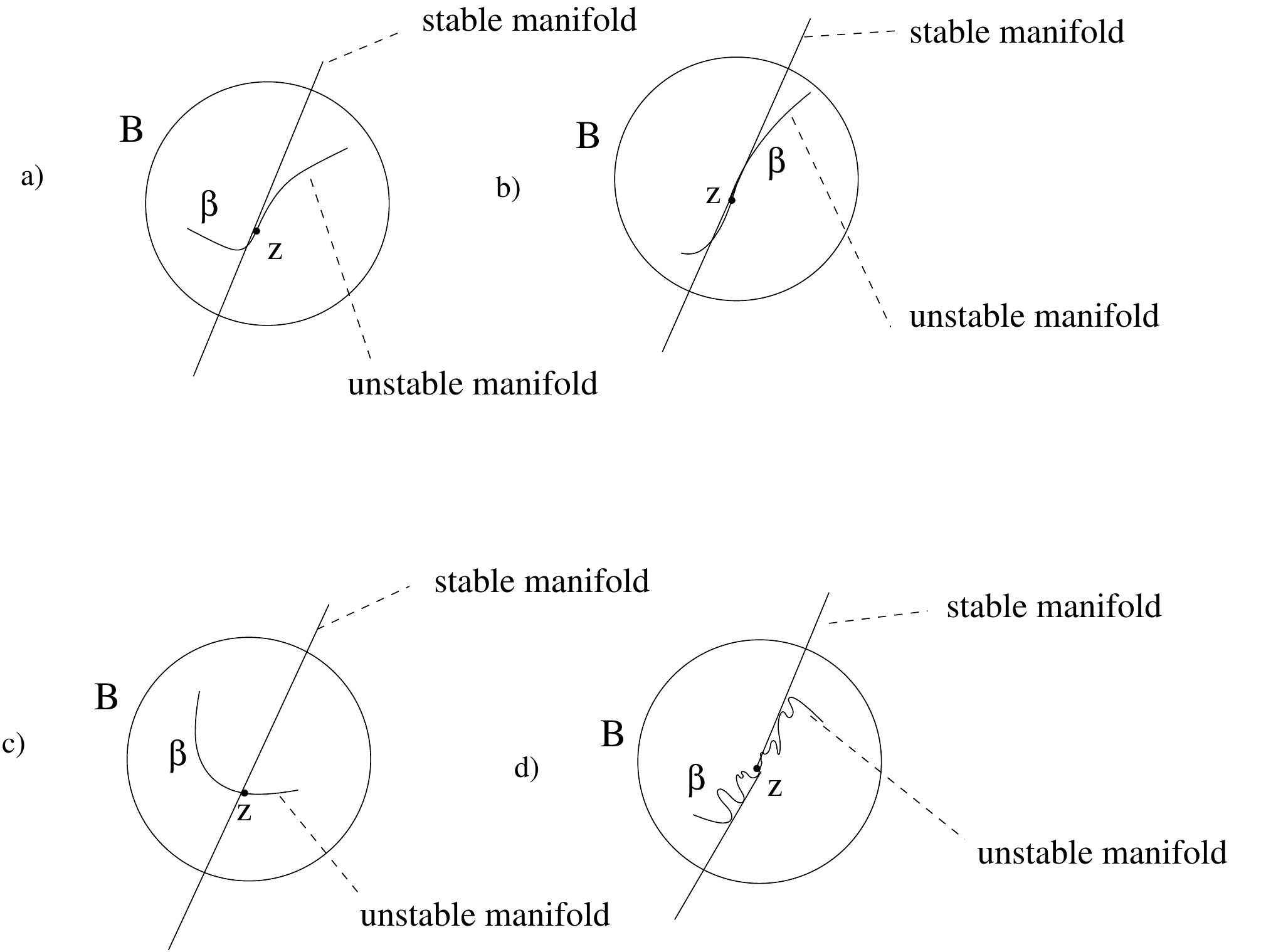}

\hfill{}

\caption{\small  4 cases of topologically transverse intersections; a) $z$ is a odd order tangency
  b) there is a segment in the intersection of the manifolds
  c) a $C^1$-transverse crossing d) $z$ is accumulated on both sides by even order tangencies.}

\label{Figure 1}

\end{figure}

The most important consequence of a topologically transverse intersection for us
is a $C^0$ $\lambda $-lemma: If $W^u(p)$ has a topologically transverse
intersection with $W^s(q),$ then $W^u(p)$ $C^0$-accumulates on $W^u(q).$

As pointed out in \cite{jlms}, the following converse of (\ref{resvida}) is
true: if $\widetilde{g}\stackrel{def.}{=}\widetilde{f}^q(\bullet )-(p,l)$
has a hyperbolic periodic saddle $\widetilde{P}\in {\rm I\negthinspace R^2}$
such that $W^u(\widetilde{P})$ has a topologically transverse intersection
with $W^s(\widetilde{P})+(a_i,b_i),$ for integer vectors $(a_i,b_i),$ $%
i=1,2,...,k,$ such that 
$$
(0,0)\in ConvexHull\{(a_1,b_1),(a_2,b_2),...,(a_k,b_k)\}, 
$$
then $(0,0)\in int(\rho (\widetilde{g}))\Leftrightarrow (p/q,l/q)\in
int(\rho (\widetilde{f})).$ This follows from the following:

\vskip0.2truecm

{\bf Lemma 0.} {\it Let $g\in Diff_0^1({\rm T^2})$ and $\widetilde{g}:{\rm I%
\negthinspace R^2\rightarrow I\negthinspace R^2}$ be a lift of $g$ which has
a hyperbolic periodic saddle point $\widetilde{P}$ such that $W^u(\widetilde{%
P})$ has a topologically transverse intersection with $W^s(\widetilde{P}%
)+(a,b),$ for some integer vector $(a,b)\neq (0,0).$ Then $\rho (\widetilde{g%
})$ contains $(0,0)$ and a rational vector parallel to $(a,b)$ with the same
orientation as $(a,b).$}

\vskip0.2truecm

In order to prove this lemma, one just have to note that if $W^u(\widetilde{P%
})$ has a topologically transverse intersection with $W^s(\widetilde{P}%
)+(a,b),$ then we can produce a topological horseshoe for $\widetilde{g}$
(see \cite{jlms}), for which a certain periodic sequence will correspond to
points with rotation vector parallel and with the same orientation as $%
(a,b). $
So, when $(p/q,l/q)\in \partial \rho (\widetilde{f})$ for some $\widetilde{f}%
\in Diff_0^2({\rm I\negthinspace R^2}),$ it may be the case that $\widetilde{%
f}^q(\bullet )-(p,l)$ has a hyperbolic periodic saddle $\widetilde{P}$ such
that $W^u(\widetilde{P})$ has a topologically transverse intersection with $%
W^s(\widetilde{P})+(a,b),$ for some integer vectors $(a,b),$ but not for all.

Moreover, if $r$ is a supporting line at $(p/q,l/q)\in \partial \rho (
\widetilde{f}),$ which means that $r$ is a straight line which contains $%
(p/q,l/q)$ and does not intersect $int(\rho (\widetilde{f})),$ and if $
\overrightarrow{v}$ is a vector orthogonal to $r,$ such that $-
\overrightarrow{v}$ points towards the rotation set, then $W^u(\widetilde{P})
$ has a topologically transverse intersection with $W^s(\widetilde{P})+(a,b),
$ for some integer vector $(a,b)$ $\Rightarrow $ $(a,b).\overrightarrow{v}%
\leq 0.$ If $\rho (\widetilde{f})$ intersects $r$ only at $(p/q,r/q),$ then $%
(a,b).\overrightarrow{v}\geq 0\Rightarrow (a,b)=(0,0).$

On the family $f_t$ and the rational $\rho =(p/q,l/q)$ which is in the
boundary of the rotation set $\rho (\widetilde{f}_t)$ at the critical
parameter $t=\overline{t},$ without loss of generality, we can assume $\rho $
to be $(0,0).$ Instead of considering $f_t$ and its lift $\widetilde{f}_t,$
we just have to consider $f_t^q$ and the lift $\widetilde{f}_t^q-(p,l).$
This is a standard procedure for this type of problem. We just have to be
careful because we will assume some hypotheses for the family $f_t,
\widetilde{f}_t$ and we must see that they also hold for $f_t^q,\widetilde{f}%
_t^q-(p,l).$ Let us look at this.

\subsubsection{Hypotheses of the main theorem}

Assume $f_t\in Diff_0^\infty ({\rm T^2})$ is a generic $C^\infty $-family of
area-preserving diffeomorphisms ($t\in ]\overline{t}-\epsilon ,\overline{t}%
+\epsilon [$ for some parameter $\overline{t}$ and $\epsilon >0$), among
other things, in the sense of Meyer \cite{meyer1970} such that:

\begin{enumerate}
\item  $\rho (\widetilde{f}_t)$ has interior for all $t\in ]\overline{t}%
-\epsilon ,\overline{t}+\epsilon [;$

\item  $(p/q,l/q)\in \partial \rho (\widetilde{f}_{\overline{t}}),$ $r$ is a
supporting line for $\rho (\widetilde{f}_{\overline{t}})$ at $(p/q,l/q),$ $
\overrightarrow{v}$ is a unitary vector orthogonal to $r,$ such that $-
\overrightarrow{v}$ points towards the rotation set;

\item  $(p/q,l/q)\in int(\rho (\widetilde{f}_t))$ for all $t\in ]\overline{t}%
,\overline{t}+\epsilon [;$

\item  The genericity in the sense of Meyer implies that if for some
parameter $t\in ]\overline{t}-\epsilon ,\overline{t}+\epsilon [,$ a $f_t$%
-periodic point has $1$ as an eigenvalue, then it is a saddle-elliptic type
of point, one which is going to give birth to a saddle and an elliptic point
when the parameter moves in one direction and it is going to disappear if
the parameter moves in the other direction. As the family is generic, for
each period there are only finitely many periodic points. Moreover, as is
explained in page 3 of the summary of \cite{dumortier}, we can assume, that
at all periodic points, in particular, at each $n$-periodic point (for all
integers $n>0$) which has $1,1$ as eigenvalues (the point is isolated among
the $n$-periodic points and must have zero topological index), a Lojasiewicz
condition is satisfied. So, as explained in section 2 of \cite{eupatrice},
the dynamics near such a point is as in figure 2.

\begin{figure}[ht!]

\hfill

\includegraphics [height=78mm]{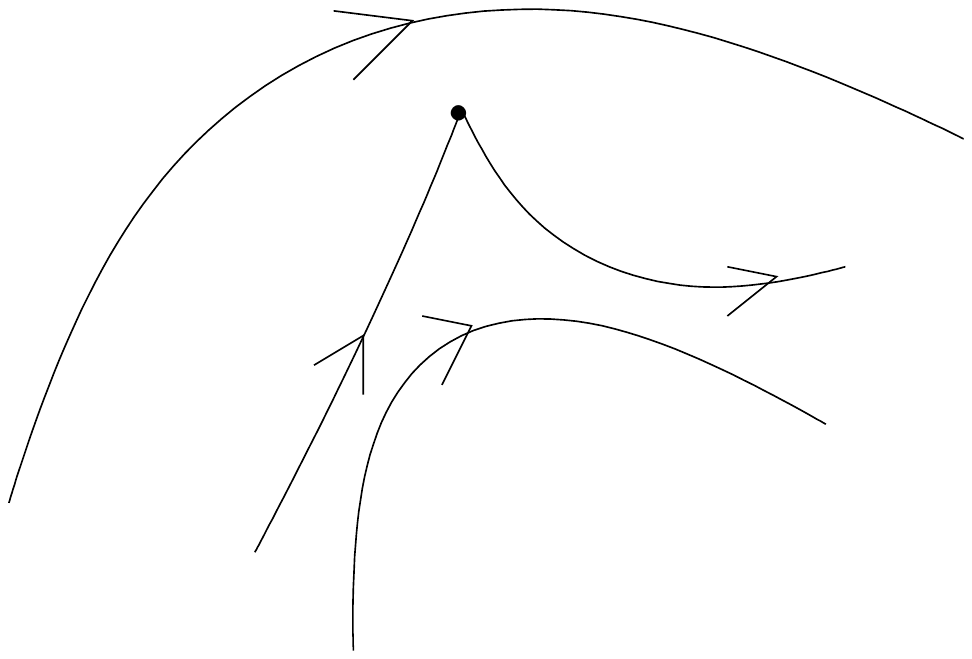}

\hfill{}

\caption{\small Dynamics in a neighborhood of the degenerate periodic points.}

\label{Figure 2}

\end{figure}

\item  Saddle-connections are a phenomena of infinity codimension (see \cite
{durmelo}). Therefore, as we are considering $C^\infty $-generic 1-parameter
families, we can also assume that for all $t\in ]\overline{t}-\epsilon ,
\overline{t}+\epsilon [,$ $f_t$ does not have connections between invariant
branches of periodic points, which can be either hyperbolic saddles or
degenerate as in figure 2. This is not strictly contained in the literature
on this subject, but a proof assuming this more general situation can be
obtained exactly in the same way as is done when only hyperbolic saddles are
considered.

\item  Moreover, as explained in section 6 of chapter II\ of \cite{nepata},
a much stronger statement holds:\ for $C^\infty $ 1-parameter generic families $%
f_t,$ if a point $z\in {\rm T^2}$ belongs to the intersection of a stable
and an unstable manifold of some $f_t$-periodic hyperbolic saddles and the
intersection is not $C^1$-transverse, then it is a quadratic tangency, that
is, it is not topologically transverse. This implies that every time an
unstable manifold has a topologically transverse intersection with a stable
manifold of some hyperbolic periodic points, this intersection is actually $%
C^1$-transverse. And when a tangency appears, it unfolds generically with
the parameter (with positive speed, see remark 6.2 of \cite{nepata}). This
means that if for some hyperbolic periodic saddles $q_t$ and $p_t$ such
that $W^u(q_{t^{\prime }})$ and $W^s(p_{t^{\prime }})$ have a quadratic
tangency at a point $z_{t^{\prime }}\in {\rm T^2}$ for some parameter $%
t^{\prime },$ then for $t$ close to $t^{\prime },$ in suitable coordinates
near $z_{t^{\prime }}=(0,0),$ we can write $W^u(q_t)=(x,f(x)+(t-t^{\prime }))
$ and $W^s(p_t)=(x,0),$ where $f$ is a $C^\infty $ function defined in a
neighborhood of $0$ such that $f(0)=0,$ $f^{\prime }(0)=0$ and $f^{\prime
\prime }(0)>0.$ This implies that if the parameter varies in a neighborhood
of the tangency parameter, to one side a $C^1$-transverse intersection is
created and to the other side, the intersection disappears.

This is stated for families of general diffeomorphisms in \cite{nepata}, but
the same result holds for families of area-preserving diffeomorphisms. It is
not hard to see that in order to avoid degenerate tangencies (of order $\geq
3$), the kind of perturbations needed can be performed preserving area.
\end{enumerate}

If $f_t,\widetilde{f}_t$ satisfy the above hypotheses, then $f_t^q,
\widetilde{f}_t^q-(p,l)$ clearly satisfy the same set of hypotheses with
respect to $\rho =(0,0).$ Note that the supporting line at $(0,0)$ for $\rho
(\widetilde{f}_{\overline{t}}^q-(p,l))$ is parallel to $r,$ the supporting line for $\rho
(\widetilde{f}_{\overline{t}})$ at $(p/q,l/q).$ So there is no restriction
in assuming $(p/q,l/q)$ to be $(0,0).$

\subsubsection{Statement of the main theorem}

\begin{theorem}
\label{main1}{\bf.} Under the 6 hypotheses assumed in 1.3.3 for the
family $f_t$ with $(p/q,l/q)=(0,0),$ the following holds: $\widetilde{f}_{
\overline{t}}$ has a hyperbolic fixed saddle $\widetilde{P}_{\overline{t}}$
such that $W^u(\widetilde{P}_{\overline{t}})$ has a topologically transverse
(and therefore a $C^1$-transverse) intersection with $W^s(\widetilde{P}_{
\overline{t}})+(a,b)$ for some integer vector $(a,b)\neq (0,0).$ And there
exists $K(f)>0$ such that for any $(c,d)\in{\rm Z\negthinspace
\negthinspace Z^2}$ for
which $(c,d).\overrightarrow{v}>K(f),$ if $\widetilde{P}_t$ is the
continuation of $\widetilde{P}_{\overline{t}}$ for $t>\overline{t},$ then
there exists a sequence $t_i>\overline{t},$ $t_i\stackrel{i\rightarrow
\infty }{\rightarrow }\overline{t}$ such that $W^u(\widetilde{P}_{t_i})$ has
a quadratic tangency with $W^s(\widetilde{P}_{t_i})+(c^{\prime },d^{\prime
}),$ for some $(c^{\prime },d^{\prime })\in {\rm Z\negthinspace
\negthinspace Z^2}$ which satisfies $\left| (c^{\prime }-c,d^{\prime }-d).
\overrightarrow{v}\right| \leq K(f)/4,$ and the tangency generically
unfolds for $t>t_i.$ The vector $(c^{\prime },d^{\prime })$ is within a
bounded distance from $(c,d)$ in the direction of $\overrightarrow{v},$ but
may be far in the direction of $\overrightarrow{v^{\perp }}.$ These
tangencies, are heteroclinic intersections for $\widetilde{f}_{t_i}$ which
could not exist at $t=\overline{t}.$ Finally, we point out that for all $t>
\overline{t}$ and for all integer vectors $(a,b),$ $W^u(\widetilde{P}_t)$
has a $C^1$ transverse intersection with $W^s(\widetilde{P}_t)+(a,b).$
\end{theorem}

{\bf Remarks:}

\begin{itemize}
\item  As we said, the tangencies given in the previous theorem are
precisely for some integer vectors $(c^{\prime },d^{\prime })$ for which at $%
t=\overline{t},$ they could not exist. This will become clear in the proof. 
%

\item  We were not able to produce tangencies at $t=\overline{t},$ even when 
$(c^{\prime },d^{\prime }).\overrightarrow{v}>0$ is small. In fact, an on
going work by Jager, Koropecki and Tal for a particular family suggests that
for their family, these tangencies may not exist at the bifurcating
parameter.\ Actually, it is easy to see that for an area preserving
diffeomorphism of the sphere, if it has a hyperbolic saddle fixed point $p$
and $W^u(p)$ intersects $W^s(p),$ then there is a topologically transverse
intersection between $W^u(p)$ and $W^s(p).$ So, for a generic family of such
maps, a horseshoe is not preceded by a tangency: the existence of a tangency
already implies a horseshoe. This is clearly not true out of the the area
preserving world and shows how subtle is the problem of birth of a horseshoe
in the conservative case.

\item  Intuitively, as the rotation set becomes larger, one would expect the
topological entropy to grow, at least for a tight model. In fact, in Kwapisz 
\cite{thesis}, some lower bounds for the topological entropy related to the
2-dimensional size of the rotation set are presented (it is conjectured
there that the area of the rotation set could be used, but what is actually
used is a more technical computation on the size of the rotation set). Our
main theorem says that every time the rotation set locally grows near a
rational point, then nearby maps must have tangencies, which generically
unfold as the parameter changes. And this is a phenomenon which is
associated to the growth of topological entropy, see \cite{tahzibi} and \cite
{sambarino}. More precisely, in the two previous papers it is proved that
generically, if a surface diffeomorphism $f$ has arbitrarily close
neighbors with larger topological entropy, then $f$ has a periodic saddle
point with a homoclinic tangency. Both these results were not stated in the
area preserving case, in fact they may not be true in the area preserving
case, but they indicate that whenever the topological entropy is ready to
grow, it is expected to find tangencies nearby.

\item  The unfolding of the above tangencies create generic elliptic
periodic points, see \cite{duarte}.

\item  An analytic version of the above theorem can also be proved. We have
to assume that:

\begin{enumerate}
\item  the family has no connections between separatrices of periodic points;

\item  for each period, there are only finitely many periodic points;

\item  if a periodic point has negative topological index, then it is a
hyperbolic saddle;
\end{enumerate}

These conditions are generic among $C^\infty $-$1$-parameter families, but
for analytic families I\ do not know. The tangencies obtained in this case
have finite order, but are not necessarily quadratic. To prove such a
result, first remember that all isolated periodic points for analytic
area-preserving diffeomorphisms satisfy a Lojasiewicz condition, see section
2 of \cite{eupatrice}. And moreover, if an isolated periodic point has a
characteristic curve (see \cite{dumortier} and again \cite{eupatrice},
section 2), then from the preservation of area, the dynamics in a
neighborhood of such a point is obtained, at least in a topological sense,
by gluing a finite number of saddle sectors.

Another important ingredient is the main result of \cite{patmeykoro2} quoted here
as theorem A,
which among other things, says that for an area-preserving diffeomorphism $f$
of the plane, which for every $n>0$ has only isolated $n$-periodic points,
if it has an invariant topological open disk $U$ with compact boundary,
whose prime ends rotation number is rational, then $\partial U$ contains
periodic points, all of the same period $k>0$ and the eigenvalues of ($Df^k)$
at these periodic points contained in $\partial U$ are real and positive.
So, if such a map $f$ is analytic, the $k$-periodic points in $\partial U$
satisfy a Lojasiewicz condition. If some of these periodic points, for
instance denoted $P,$ has topological index $1,$ then the
eigenvalues of ($Df^k\mid _P)$ are both equal to $1.$ Under these
conditions, the main result of \cite{jiang} implies the existence of
periodic orbits rotating around $P$ with many different velocities with
respect to an isotopy $I_t$ from the $Id$ to $f.$ And a technical result in 
\cite{patmeykoro2} says that if a periodic point belongs to $\partial U,$
then this can not happen. So, the topological indexes of all periodic points
in $\partial U$ are less or equal to zero and thus from section 2 of \cite{eupatrice},
they all have characteristic curves. Therefore, locally, all periodic points in
$\partial U$ are saddle like. They may have 2 sectors (index zero) or 4 sectors
(index -1). And from \cite{mather}, if 
$\partial U$ is bounded, connections must exist. Thus, the hypothesis that
there are no connections between separatrices of periodic points imply the
irrationality of the prime ends rotation number for all open invariant disks
whose boundaries are compact.

And finally, the last result we need is due to Churchill and Rod \cite
{churrod}. It says that, for analytic area preserving diffeomorphisms, the
existence of a topologically transverse homoclinic point for a certain
saddle, implies the existence of a $C^1$-transverse homoclinic point for
that saddle. Using these results in the appropriate places of the proof in
the next section, an analytic version of the main theorem can be obtained.
\end{itemize}

In the next section of this paper we prove our main result.

\section{Proof of the main theorem}

The proof will be divided in 2 steps.

\subsection{Step 1}

Here we prove that $\widetilde{f}_{\overline{t}}$ has a hyperbolic fixed
saddle such that its unstable manifold has a topologically transverse
intersection (therefore, $C^1$-transverse) with its stable manifold
translated by a non-zero integer vector $(a,b).$ Clearly, from lemma 0, $%
(a,b).\overrightarrow{v}\leq 0.$

First of all, note that as $(0,0)\in \partial \rho (\widetilde{f}_{\overline{%
t}})$ and for all $t>\overline{t},$ $(0,0)\in int(\rho (\widetilde{f}_t)),$ $
\widetilde{f}_{\overline{t}}$ has (finitely many) fixed points up to ${\rm Z%
\negthinspace
\negthinspace Z^2}$ translations, it can not be fixed point free. The
finiteness comes from the generic assumptions. If all these fixed points had
zero topological index, as is explained before the statement of theorem 1,
the dynamics near each of them would be as in figure 2. And in this
situation, theorem 1 of \cite{eupatrice} implies that $(0,0)\notin int(\rho (
\widetilde{f}_t))$ for any $t$ close to $\overline{t}.$

So there must be $\widetilde{f}_{\overline{t}}$-fixed points with non-zero
topological index. From the Nielsen-Lefschetz index theorem, we obtain a
fixed point for $\widetilde{f}_{\overline{t}}$ with negative index. From the
genericity of our family, the only negative index allowed is $-1$ and fixed
points with topological index equal to $-1$ are hyperbolic saddles. Assume
there are $k>0$ of such points in the fundamental domain $[0,1[^2,$ denoted $%
\{\widetilde{P}_{\overline{t}}^1,...,\widetilde{P}_{\overline{t}}^k\}.$ So,
in $[0,1[^2$ $\widetilde{f}_{\overline{t}}$ has $k$ hyperbolic fixed saddle
points and other periodic points with topological index greater or equal to
zero.

Now let us choose a rational vector in $int(\rho (\widetilde{f}_{\overline{t}%
})).$ Without loss of generality, conjugating $f$ with some adequate integer
matrix if necessary, we can suppose that this rational vector is of the form 
$(0,-1/n)$ for some $n>0.$

By some results from \cite{c1epsilon}, let $\widetilde{Q}\in {\rm I%
\negthinspace R^2}$ be a periodic hyperbolic saddle point for $\left( 
\widetilde{f}_{\overline{t}}\right) ^n+(0,1)$ such that $W^u(\widetilde{Q})$
has a topologically transverse intersection with $W^s(\widetilde{Q})+(a,b)$
for all integer vectors $(a,b).$ In this case, $\overline{W^u(\widetilde{Q})}%
=\overline{W^s(\widetilde{Q})}=R.I.(\widetilde{f}_{\overline{t}})=Region$ $%
of $ $instability$ $of$ $\widetilde{f}_{\overline{t}},$ a $\widetilde{f}_{
\overline{t}}$-invariant equivariant set such that, if $\widetilde{D}$ is a
connected component of its complement, then $\widetilde{D}$ is a connected
component of the lift of a $f_{\overline{t}}$-periodic open disk in the
torus and for every $f_{\overline{t}}$-periodic open disk $D\subset {\rm T^2,%
}$ $p^{-1}(D)\subset \left( R.I.(\widetilde{f}_{\overline{t}})\right) ^c.$
Also, for any rational vector $\left( p/q,l/q\right) \in int(\rho ( 
\widetilde{f}_{\overline{t}})),$ there exists a hyperbolic periodic saddle
point for $\left( \widetilde{f}_{\overline{t}}\right) ^q-(p,l),$ such that
its unstable manifold also has topologically transverse intersections with
all integer translates of its stable manifold and so, the closure of its
stable manifold is equal the closure of its unstable manifold and they are
both equal $R.I.(\widetilde{f}_{\overline{t}}).$ As we said, these results were proved
in \cite{c1epsilon} and similar statements hold for homeomorphisms \cite{talkoro}.

If for some $1\leq i^{*}\leq k,$ $W^u(\widetilde{P}_{\overline{%
t}}^{i^{*}})$ and $W^s(\widetilde{P}_{\overline{t}}^{i^{*}})$ are both
unbounded subsets of the plane, then it can be proved that $W^u(\widetilde{P}%
_{\overline{t}}^{i^{*}})$ must have a topologically transverse intersection
with $W^s(\widetilde{Q})$ and $W^s(\widetilde{P}_{\overline{t}}^{i^{*}})$
must have a topologically transverse intersection with $W^u(\widetilde{Q}).$
But, as the rotation vector of $\widetilde{Q}$ is not zero, this
gives what we want in step 1. More precisely, the following fact
holds:

\begin{description}
\item[Fact]  :{\it \ If }$W^u(\widetilde{P}_{\overline{t}}^{i^{*}})$ and $%
W^s(\widetilde{P}_{\overline{t}}^{i^{*}})$ are both unbounded subsets of the
plane, then $W^u(\widetilde{P}_{\overline{t}}^{i^{*}})$ has a topologically
transverse intersection with $W^s(\widetilde{P}_{\overline{t}%
}^{i^{*}})-(0,1).$

{\it Proof:}
\end{description}

As $W^s(\widetilde{Q})$ has a topologically transverse intersection with $%
W^u(\widetilde{Q})+(a,b)$ for all integer vectors $(a,b),$ this implies that
if $W^u(\widetilde{P}_{\overline{t}}^{i^{*}})$ is unbounded, then $W^u( 
\widetilde{P}_{\overline{t}}^{i^{*}})$ has a topologically transverse
intersection with $W^s(\widetilde{Q}).$ This follows from the following
idea: There is a compact arc $\lambda _u$ in $W^u(\widetilde{Q})$ that
contains $\widetilde{Q}$ and a compact arc $\lambda _s$ in $W^s(\widetilde{Q}%
)$ that also contains $\widetilde{Q},$ such that $\lambda _u$ has
topologically transverse intersections with $\lambda _s+(0,1)$ and $\lambda
_s+(1,0).$ This implies that the connected components of the complement of 
$$
\cup _{(a,b)\in {\rm Z\negthinspace
\negthinspace Z^2}}\lambda _u\cup \lambda _s+(a,b) 
$$
are all open topological disks, with diameter uniformly bounded from above.
So, if $W^u(\widetilde{P}_{\overline{t}}^{i^{*}})$ is unbounded, it must
have a topologically transverse intersection with some translate of $\lambda
_s.$ As $W^s(\widetilde{Q})$ $C^0$-accumulates on all its integer
translates, we finally get that $W^u(\widetilde{P}_{\overline{t}}^{i^{*}})$
has a topologically transverse intersection with $W^s(\widetilde{Q}).$

So, as for some integer $m>0,$%
$$
\left( \widetilde{f}_{\overline{t}}\right) ^{m.n}(\widetilde{Q})=\widetilde{Q%
}-(0,m), 
$$
$W^u(\widetilde{P}_{\overline{t}}^{i^{*}})$ $C^0$-accumulates on compact
pieces of $W^u(\widetilde{Q})-(0,k_j)$ for a certain sequence $%
k_j\rightarrow \infty ,$ that is, given a compact arc $\widetilde{\theta }$
contained in $W^u(\widetilde{Q}),$ there exists a sequence $k_j\rightarrow
\infty $ such that for some arcs $\widetilde{\theta }_j\subset W^u(\widetilde{%
P}_{\overline{t}}^{i^{*}}),$ $\widetilde{\theta }_j+(0,k_j)\rightarrow 
\widetilde{\theta }$ in the Hausdorff topology as $j\rightarrow \infty .$ An
analogous argument implies that if $W^s(\widetilde{P}_{\overline{t}%
}^{i^{*}}) $ is unbounded, then $W^s(\widetilde{P}_{\overline{t}}^{i^{*}})$
has a topologically transverse intersection with $W^u(\widetilde{Q}).$ So if
we choose a compact arc $\kappa _u$ contained in $W^u(\widetilde{Q})$ which
has a topologically transverse intersection with $W^s(\widetilde{P}_{
\overline{t}}^{i^{*}}),$ we get that $W^u(\widetilde{P}_{\overline{t}%
}^{i^{*}})$ accumulates on $\kappa _u-(0,k_j)$ and thus it has a
topologically transverse intersection with $W^s(\widetilde{P}_{\overline{t}%
}^{i^{*}})-(0,k_j)$ for some $k_j>0$ sufficiently large$.$ As we pointed out
after the definition of topologically transverse intersections, all
the above follows from a $C^0$-version of the $\lambda $-lemma that holds
for topologically transverse intersections. Now consider a compact subarc
of a branch of $W^u(\widetilde{P}_{\overline{t}}^{i^{*}}),$ denoted $\alpha
^u,$ starting at $\widetilde{P}_{\overline{t}}^{i^{*}}$ and a compact subarc
of a branch of $W^s(\widetilde{P}_{\overline{t}}^{i^{*}}),$ denoted $\alpha
^s,$ starting at $\widetilde{P}_{\overline{t}}^{i^{*}}$ such that 
\begin{equation}
\label{toptrans1} 
\begin{array}{c}
\alpha ^u 
\text{ has a topologically transverse intersection with} \\ \alpha
^s-(0,k_j) \text{ for some }k_j>0. 
\end{array}
\end{equation}
Using Brouwer's lemma on translation arcs exactly as in lemma 24 of \cite
{c1epsilon}, we get that either $\alpha ^u$ has an intersection with $\alpha
^s-(0,1)$ or $\alpha ^u-(0,1)$ has an intersection with $\alpha ^s.$ If $%
\alpha ^u$ had only non topologically transverse intersections with $\alpha
^s-(0,1)$ and $\alpha ^u-(0,1)$ had only non topologically transverse
intersections with $\alpha ^s,$ then we could $C^0$-perturb $\alpha ^u$ and $%
\alpha ^s$ in an arbitrarily small way such that $(\alpha _{per}^u\cup
\alpha _{per}^s)\cap ((\alpha _{per}^u\cup \alpha _{per}^s)-(0,1))=\emptyset 
$ and $\alpha _{per}^u\cap (\alpha _{per}^s-(0,k_j))\neq \emptyset $
(because of the topologically transverse assumption (\ref{toptrans1})). But
this contradicts Brouwer 's lemma \cite{brown}. So, either $\alpha ^u$ has a
topologically transverse intersection with $\alpha ^s-(0,1)$ or $\alpha
^u-(0,1)$ has a topologically transverse intersection with $\alpha ^s.$ The
first possibility is what we want and the second is not possible because
lemma 0 would imply that $\rho (\widetilde{f}_{\overline{t}})$ contains a
point of the form $(0,a)$ for some $a>0.$ As $(0,-1/n)\in int(\rho
( \widetilde{f}_{\overline{t}}))$ these two facts contradict the assumption
that $(0,0)\in \partial \rho (\widetilde{f}_{\overline{t}}).$ $\Box $

\vskip 0.2truecm

In order to conclude this step, we need the following lemma:

\begin{lemma}
\label{fimstep1}. There exists $1\leq i^{*}\leq k$ such that for any choice
of $\lambda _u^{i^{*}}$ and $\lambda _s^{i^{*}},$ one unstable and one
stable branch at $\widetilde{P}_{\overline{t}}^{i^{*}},$ they are both
unbounded.
\end{lemma}

{\it Proof: }

For each $1\leq i\leq k,$ as $\widetilde{P}_{\overline{t}}^i$ has
topological index equal $-1,$ the two stable and the two unstable branches
are each, $\widetilde{f}_{\overline{t}}$-invariant. Fixed some unstable
branch $\lambda _u^i$ and some stable $\lambda _s^i,$ let 
$$
K_u^i=\overline{\lambda _u^i}\text{ and }K_s^i=\overline{\lambda _s^i}.\text{
} 
$$
Both are connected $\widetilde{f}_{\overline{t}}$-invariant sets. Let $K^i$
be equal to either $K_u^i$ or $K_s^i$ and assume it is bounded. Without loss
of generality, suppose that $K^i=K_u^i.$ First, we collect some properties
about $K^i.$

\begin{itemize}
\item  If $K^i$ intersects a connected component $\widetilde{D}$ of the
complement of $R.I.(\widetilde{f}_{\overline{t}}),$ then from lemma 6.1 of 
\cite{frankslecal}, $\lambda _u^i\backslash \widetilde{P}_{\overline{t}}^i$
is contained in $\widetilde{D},$ which then, is a $\widetilde{f}_{\overline{t%
}}$-invariant bounded open disk (see \cite{c1epsilon} and \cite{talkoro}).
As the family of diffeomorphisms considered is generic (in particular, it
does not have connections between stable and unstable separatrices of
periodic points), the rotation number of the prime ends compactification of $
\widetilde{D},$ denoted $\beta ,$ must be irrational by theorem A. 
So, if $\widetilde{P}_{\overline{t}}^i\in \partial $$
\widetilde{D},$ as $\lambda _u^i\backslash \widetilde{P}_{\overline{t}%
}^i\subset \widetilde{D},$ this would be a contradiction with the
irrationality of $\beta ,$ because $\widetilde{f}_{\overline{t}}(\lambda
_u^i)=\lambda _u^i.$ Thus, $\widetilde{P}_{\overline{t}}^i$ is contained in $
\widetilde{D}$ and the topological index of $\widetilde{f}_{\overline{t}}$
with respect to $\widetilde{D}$ is $+1$ (because $\beta $ is irrational,
therefore not zero). By topological index of $\widetilde{f}_{\overline{t}}$
with respect to $\widetilde{D}$ we mean the sum of the indexes at all the $
\widetilde{f}_{\overline{t}}$-fixed points contained in $\widetilde{D}.$ This
information will be used in the end of the proof.
\end{itemize}

%
%

\begin{itemize}
\item  Suppose now that $K^i=K_u^i$ is contained in $R.I.(\widetilde{f}_{
    \overline{t}}).$ It is not possible that $K_u^i\cap \lambda _s^i=\widetilde{P}_{\overline{t}}^i$
because it would imply that the connected component $M$ of $(K^i)^c$ which
contains $\lambda _s^i\backslash \widetilde{P}_{\overline{t}}^i$ has
rational prime ends rotation number and as we already said, this does not
happen under our generic conditions. So, from lemma 2 of Oliveira \cite
{oliveira}, we get that either $K_u^i\supset \lambda _s^i$ or $\lambda _u^i$
intersects $\lambda _s^i.$ If $\lambda _u^i$ intersects $\lambda _s^i,$ then 
$\lambda _s^i$ intersects $R.I.(\widetilde{f}_{\overline{t}}),$ so it is
contained in $R.I.(\widetilde{f}_{\overline{t}})$ and we can find a Jordan
curve $\tau $ contained in $\lambda _u^i\cup \lambda _s^i,$ $\widetilde{P}_{
\overline{t}}^i\in \tau .$ Theorem B implies that there is
no periodic point in the boundary of a connected component of the complement
of $R.I.(\widetilde{f}_{\overline{t}}),$ because such a component is $f_{
\overline{t}}$-periodic when projected to the torus and it has irrational
prime ends rotation number by theorem A. So, interior($%
\tau $) intersects $R.I.(\widetilde{f}_{\overline{t}}),$ and thus interior($%
\tau $) intersects both $W^u(\widetilde{Q})$ and $W^s(\widetilde{Q}),$ something
that contradicts the assumption that $K_u^i$ is bounded because $\lambda _u^i$
must intersect $W^s(\widetilde{Q}),$ which implies that it is unbounded.

And if $K_u^i\supset \lambda _s^i,$ then from our assumption that $K_u^i$ is
bounded, we get that $K_s^i=\overline{\lambda _s^i}$ is also bounded.
Arguing as above, we obtain that $K_s^i\cap \lambda _u^i\neq \widetilde{P}_{\overline{t}}^i.$

Otherwise, if $M^{*}$ is the connected component of $(K_s^i)^c$ that
contains $\lambda _u^i\backslash \widetilde{P}_{\overline{t}}^i,$ as $
\widetilde{f}_{\overline{t}}(\lambda _u^i)=\lambda _u^i,$ the rotation
number of the prime ends compactification of $M^{*}$ would be rational. As
this implies connections between separatrices of periodic points, which do
not exist under our hypotheses, $K_s^i$ must intersect
$\lambda _u^i\backslash \widetilde{P}_{\overline{t}}^i.$ And
again, from lemma 2 of Oliveira \cite{oliveira}, $K_s^i\supset \lambda _u^i$
or $\lambda _u^i$ intersects $\lambda _s^i$ and we are done.

Thus, the almost final situation we have to deal is when $K_u^i\supset
\lambda _s^i$ and $K_s^i\supset \lambda _u^i.$ But it is contained in the
proof of the main theorem of \cite{oliveira} that these relations imply that 
$\lambda _u^i$ intersects $\lambda _s^i.$ And as we explained above, this is
a contradiction with the assumption that $K_u^i$ is bounded. If $K^i=K_s^i$
and $K_s^i\subset R.I.(\widetilde{f}_{\overline{t}}),$ an analogous argument
could be applied in order to arrive at similar contradictions. 
\end{itemize}
Thus, if $K^i$ is bounded, it must be contained in the complement
of $R.I.(\widetilde{f}_{\overline{t}}).$ 
In order to conclude the proof, we are left to consider the case when
for every $1\leq i\leq k,$ we can choose $K^i$ either equal to $K_u^i$ or $%
K_s^i,$ such that it is bounded and contained in a connected component $
\widetilde{D}_i$ of the complement of $R.I.(\widetilde{f}_{\overline{t}}).$
As we already obtained, the topological index of $\widetilde{f}_{\overline{t}%
}$ restricted to $\widetilde{D}_i$ is $+1.$ This clearly contradicts the
Nielsen-Lefschetz index formula because the sum of the indices of $f_{
\overline{t}}$ at its fixed points which have $(0,0)$ rotation vector would
be positive. So, for some $1\leq i^{*}\leq k,$ both unstable and both stable
branches at $\widetilde{P}_{\overline{t}}^{i^{*}}$ are unbounded. $\Box $

\vskip 0.2truecm

This concludes step 1.

\subsection{Step 2}

From the previous step we know that there exists a $\widetilde{f}_{\overline{%
t}}$-fixed point denoted $\widetilde{P}_{\overline{t}}$ such that $W^u( 
\widetilde{P}_{\overline{t}})$ has a topologically transverse, and therefore
a $C^1$-transverse intersection with $W^s(\widetilde{P}_{\overline{t}%
})-(0,1).$

So, there exists a compact connected piece of a branch of $W^u(\widetilde{P}%
_{\overline{t}}),$ starting at $\widetilde{P}_{\overline{t}},$ denoted $%
\lambda _u^{\overline{t}}$ and a compact connected piece of a branch of $%
W^s( \widetilde{P}_{\overline{t}}),$ starting at $\widetilde{P}_{\overline{t}%
},$ denoted $\lambda _s^{\overline{t}},$ such that $\lambda _u^{\overline{t}%
}\cup \left( \lambda _s^{\overline{t}}-(0,1)\right) $ contains a continuous
curve connecting $\widetilde{P}_{\overline{t}}$ to $\widetilde{P}_{\overline{%
t}}-(0,1).$ The end point of $\lambda _u^{\overline{t}},$ denoted $w,$
belongs to $\lambda _s^{\overline{t}}-(0,1)$ and it is a $C^1$-transverse
heteroclinic point. The main consequence of the above is the following:

\begin{description}
\item[Fact]  : The curve $\widetilde{\gamma }_V^{\overline{t}}$ connecting $
\widetilde{P}_{\overline{t}}$ to $\widetilde{P}_{\overline{t}}-(0,1)$ contained in 
$\lambda _u^{\overline{t}}\cup \left( \lambda _s^{\overline{t}%
}-(0,1)\right) $ projects to a (not necessarily simple) closed curve in the
torus, homotopic to $(0,-1)$ and it has a continuous continuation for $t\geq 
\overline{t}$ suff. small. That is, for $t-\overline{t}\geq 0$ suff. small,
there exists a curve $\widetilde{\gamma }_V^t$ connecting $\widetilde{P}_t$
to $\widetilde{P}_t-(0,1)$ made by a piece of an unstable branch of $W^u(
\widetilde{P}_t)$ and a piece of a stable branch of $W^s(\widetilde{P}%
_t)-(0,1)$ such that $t\rightarrow $ $\widetilde{\gamma }_V^t$ is continuous
for $t-\overline{t}\geq 0$ suff. small.
\end{description}

{\it Proof:}

Immediate from the fact that $w$ is a $C^1$-transversal heteroclinic point
which has a continuous continuation for all diffeomorphisms $C^2$-close to $
\widetilde{f}_{\overline{t}}.$ $\Box $ %

\vskip 0.2truecm

As $(0,-1/n)$ is contained in $int(\rho (\widetilde{f}_{\overline{t}})),$
there are rational points in $int(\rho (\widetilde{f}_{\overline{t}}))$ with
positive horizontal coordinates. Thus, if we define 
$$
\Gamma _{V,a}^{\overline{t}}\stackrel{def.}{=}...\cup \widetilde{\gamma }_V^{
\overline{t}}+(a,2)\cup \widetilde{\gamma }_V^{\overline{t}}+(a,1)\cup 
\widetilde{\gamma }_V^{\overline{t}}+(a,0)\cup \widetilde{\gamma }_V^{
\overline{t}}+(a,-1)\cup ...,\text{ for }a\in {\rm Z\negthinspace
\negthinspace Z}, 
$$
we get that for all sufficiently large integer $n>0,$ $(\widetilde{f}_{
\overline{t}})^n(\Gamma _{V,0}^{\overline{t}})$ intersects $\Gamma _{V,1}^{
\overline{t}}$ transversely and moreover, we obtain a curve $\widetilde{%
\gamma }_H^{\overline{t}}$ connecting $\widetilde{P}_{\overline{t}}$ to $
\widetilde{P}_{\overline{t}}+(1,0)$ of the following form: it starts at $
\widetilde{P}_{\overline{t}},$ goes through the branch of $W^u(\widetilde{P}%
_{\overline{t}})$ which contains $\lambda _u^{\overline{t}}$ until it hits
in a topologically transverse way (so in a $C^1$-transverse way) $\lambda
_s^{\overline{t}}+(1,b)$ for some integer $b.$ If $b<0,$ we add to this
curve the following one: 
\begin{equation}
\label{bmenorz}\widetilde{\gamma }_V^{\overline{t}}+(1,0)\cup \widetilde{%
\gamma }_V^{\overline{t}}+(1,-1)\cup ...\cup \widetilde{\gamma }_V^{
\overline{t}}+(1,b+1) 
\end{equation}
and if $b\geq 0,$ we add the following curve: 
\begin{equation}
\label{bmaiorz}\widetilde{\gamma }_V^{\overline{t}}+(1,1)\cup \widetilde{%
\gamma }_V^{\overline{t}}+(1,2)\cup ...\cup \widetilde{\gamma }_V^{\overline{%
t}}+(1,b)\cup \lambda _s^{\overline{t}}+(1,b) 
\end{equation}

\begin{figure}[ht!]

\hfill

\includegraphics [height=78mm]{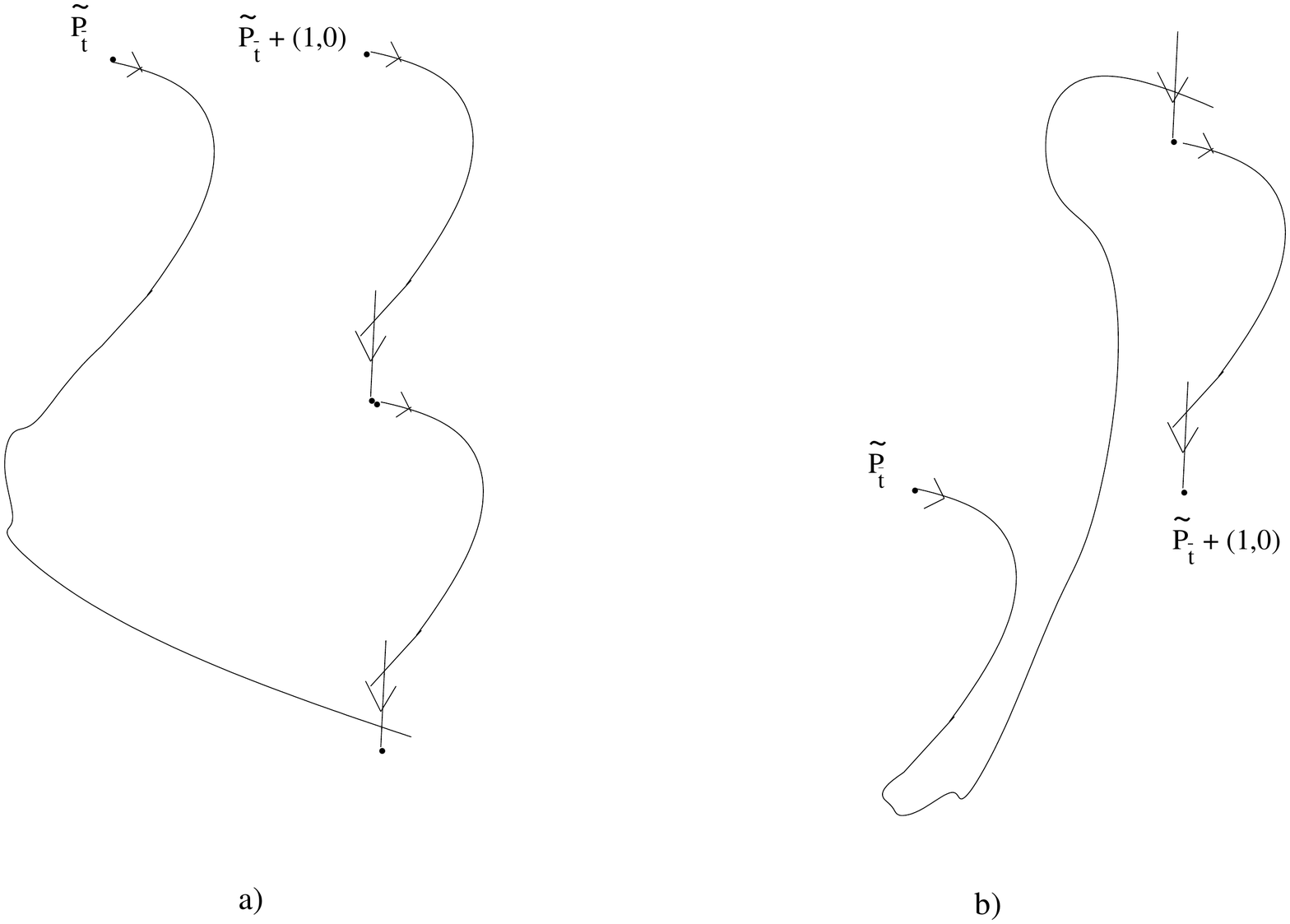}

\hfill{}

\caption{\small  How to construct $\widetilde{\gamma }_H^t.$}

\label{Figure 3}

\end{figure}

In both cases above, we omit a small piece of $\lambda _s^{\overline{t}%
}+(1,b)$ in order to get a proper curve connecting $\widetilde{P}_{\overline{%
    t}}$ and $\widetilde{P}_{\overline{t}}+(1,0).$ This follows from the fact that
when we consider iterates $(\widetilde{f}_{\overline{t}})^n(\Gamma _{V,0}^{\overline{t}})$
for some large $n>0,$ the arcs contained in stable manifolds are shrinking and arcs contained
in unstable manifolds are getting bigger. As there are orbits moving to the right under positive
iterates of $\widetilde{f}_{\overline{t}},$ the above holds.
So, $\widetilde{\gamma }_H^{
\overline{t}}$ is a continuous curve whose endpoints are $\widetilde{P}_{
\overline{t}}$ and $\widetilde{P}_{\overline{t}}+(1,0)$ and it is made of a
connected piece of an unstable branch of $W^u(\widetilde{P}_{\overline{t}}),$
added to either one of the two vertical curves above, (\ref{bmenorz}) or (%
\ref{bmaiorz}), a small piece of $\lambda _s^{\overline{t}}+(1,b)$ deleted.
As the intersection between the branch of $W^u(\widetilde{P}_{\overline{t}})$
which contains $\lambda _u^{\overline{t}}$ with $\lambda _s^{\overline{t}%
}+(1,b)$ is $C^1$-transverse and $t\rightarrow $ $\widetilde{\gamma }_V^t$
is continuous for $t-\overline{t}\geq 0$ suff. small, we get that $%
t\rightarrow $ $\widetilde{\gamma }_H^t$ is also continuous for $t-\overline{%
t}\geq 0$ suff. small. See figure 3 a) and b) for representations of these
two possibilities.

Remember that $r$ is the supporting line at $(0,0)\in \partial \rho (
\widetilde{f}_{\overline{t}})$ and $\overrightarrow{v}$ is an unitary vector
orthogonal to $r$ such that $-\overrightarrow{v}$ points towards the
rotation set.

Now we state a more general version of lemma 6 of \cite{jlms}. It does not
appear like this in that paper, but the proof presented there also proves
this more abstract version.

\vskip 0.2truecm

{\bf Lemma 6 of [2].} {\it Let $K_H$ and $K_V$ be two continua in the plane,
such that $K_H$ contains $(0,0)$ and $(1,0)$ and $K_V$ contains $(0,0)$ and $%
(0,1).$ For every vector $\overrightarrow{w},$ it is possible to construct a
connected closed set $M_{\overrightarrow{w}}$ which is equal to the union of
well chosen integer translates of $K_H$ and $K_V$ such that:

\begin{enumerate}
\item  $M_{\overrightarrow{w}}$ intersects every straight line orthogonal to 
$\overrightarrow{w};$

\item  $M_{\overrightarrow{w}}$ is bounded in the direction orthogonal to $
\overrightarrow{w},$ that is, $M_{\overrightarrow{w}}$ is contained between
two straight lines $r_{M_{\overrightarrow{w}}}$ and $s_{M_{\overrightarrow{w}%
}},$ both parallel to $\overrightarrow{w},$ and the distance between these
lines is less then $3+2.\max \{diameter(K_H),diameter(K_V)\}.$ So, in
particular $(M_{\overrightarrow{w}})^c$ has at least two unbounded connected
components, one containing $r_{M_{\overrightarrow{w}}},$ denoted $U_{r}(M_{
\overrightarrow{w}})$ and the other containing $s_{M_{\overrightarrow{w}}},$
denoted $U_{s}(M_{\overrightarrow{w}});$
\end{enumerate}
}

\vskip 0.2truecm

So, applying this lemma to the setting of this paper, it gives us a path
connected closed set $\theta _{\overrightarrow{v^{\perp }}}^{\overline{t}%
}\subset {\rm I}\negthinspace {\rm R^2}$ which is obtained as the union of certain integer
translates of $\widetilde{\gamma }_V^{\overline{t}}$ or $\widetilde{\gamma }%
_H^{\overline{t}}$ in a way that: 

\begin{enumerate}
\item  $\theta _{\overrightarrow{v^{\perp }}}^{\overline{t}}$ intersects
every straight line parallel to $\overrightarrow{v};$

\item  $\theta _{\overrightarrow{v^{\perp }}}^{\overline{t}}$ is bounded in
the direction of $\overrightarrow{v},$ that is, $\theta _{\overrightarrow{%
v^{\perp }}}^{\overline{t}}$ is contained between two straight lines $l_{-}$
and $l_{+},$ both parallel to $\overrightarrow{v^{\perp }},$ and the
distance between these lines is less then $3+2.\max \{diameter(\widetilde{%
\gamma }_V^{\overline{t}}),diameter(\widetilde{\gamma }_H^{\overline{t}})\}.$
And $(\theta _{\overrightarrow{v^{\perp }}}^{\overline{t}})^c$ has at least
two unbounded connected components, one containing $l_{-},$ denoted $U_{\_}$
and the other containing $l_{+},$ denoted $U_{+};$
\end{enumerate}

Assume that $l_{+}$ and $U_{+}$ were chosen in a way that if $(c,d)$ is an
integer vector such that $\theta _{\overrightarrow{v^{\perp }}}^{\overline{t}%
}+(c,d)$ belongs to $U_{+},$ then  
\begin{equation}
\label{escalarpos}(c,d).\overrightarrow{v}>0.
\end{equation}
It is not hard to see that for integer vectors $(c,d)$ for which $\theta _{
\overrightarrow{v^{\perp }}}^{\overline{t}}+(c,d)$ belongs either to $U_{+}$
or $U_{-},$ an inequality like (\ref{escalarpos}) needs to hold. 

For this, note that if $\overrightarrow{v}$ is a rational direction for
which $(c,d).\overrightarrow{v}=0,$ then from the way $\theta _{
\overrightarrow{v^{\perp }}}^{\overline{t}}$ is constructed, we get that $%
\theta _{\overrightarrow{v^{\perp }}}^{\overline{t}}+(c,d)$ intersects $%
\theta _{\overrightarrow{v^{\perp }}}^{\overline{t}},$ so $\theta _{
\overrightarrow{v^{\perp }}}^{\overline{t}}+(c,d)$ does not belong to $%
(\theta _{\overrightarrow{v^{\perp }}}^{\overline{t}})^c.$ And in case $
\overrightarrow{v}$ is an irrational direction, $(c,d).\overrightarrow{v}%
\neq 0.$ In particular, an integer vector $(c,d)$ satisfying the inequality
in (\ref{escalarpos}) has the property that any positive multiple of it does
not belong to $\rho (\widetilde{f}_{\overline{t}}).$ This will be important
soon.

Moreover, if $(c,d).\overrightarrow{v}>3+2.\max \{diameter(\widetilde{\gamma 
}_V^{\overline{t}}),diameter(\widetilde{\gamma }_H^{\overline{t}})\},$ then $%
\theta _{\overrightarrow{v^{\perp }}}^{\overline{t}}+(c,d)$ belongs to $%
U_{+}.$ Now we are ready to finish the proof of the main theorem.

As both $t\rightarrow $ $\widetilde{\gamma }_V^t$ and $t\rightarrow $ $
\widetilde{\gamma }_H^t$ are continuous for $t-\overline{t}\geq 0$ suff.
small, the same holds for $t\rightarrow \theta _{\overrightarrow{v^{\perp }}%
}^t.$ Given any integer vector $(c,d)$ such that 
$$
(c,d).\overrightarrow{v}>K(f)\stackrel{def.}{=}4.(3+2.\max \{diameter( 
\widetilde{\gamma }_V^{\overline{t}}),diameter(\widetilde{\gamma }_H^{
\overline{t}})\})+10, 
$$
if a $t^{*}>\overline{t},$ sufficiently close to $\overline{t}$ is fixed, we
can assume that 
$$
\theta _{\overrightarrow{v^{\perp }}}^t+(c,d)\cap \theta _{\overrightarrow{%
v^{\perp }}}^t=\emptyset \text{ for all }t\in [\overline{t},t^{*}]. 
$$
And as $(0,0)\in int(\rho (\widetilde{f}_{t^{*}})),$ there exists an integer 
$N>0$ such that 
$$
(\widetilde{f}_{t^{*}})^N(\theta _{\overrightarrow{v^{\perp }}}^{t^{*}}) 
\text{ has a topologically transverse intersection with }\theta _{
\overrightarrow{v^{\perp }}}^{t^{*}}+(c,d).\text{ } 
$$

But as $(0,0)\notin int(\rho (\widetilde{f}_{\overline{t}})),$ and $(c,d).
\overrightarrow{v}$ is sufficiently large, we get that 
$$
(\widetilde{f}_{\overline{t}})^N(\theta _{\overrightarrow{v^{\perp }}}^{
\overline{t}})\cap \theta _{\overrightarrow{v^{\perp }}}^{\overline{t}%
}+(c,d)=\emptyset .\text{ } 
$$
The previous property follows from the existence of $(a,b)\in {\rm Z%
\negthinspace
\negthinspace Z^2}$ such that 
$$
\begin{array}{c}
\theta _{
\overrightarrow{v^{\perp }}}^{\overline{t}}+(a,b)\text{ is contained between 
}\theta _{\overrightarrow{v^{\perp }}}^{\overline{t}}\text{ and }\theta _{
\overrightarrow{v^{\perp }}}^{\overline{t}}+(c,d)\\ 
\text{and} \\ (a,b).\overrightarrow{v}>2.(3+2.\max \{diameter(
\widetilde{\gamma }_V^{\overline{t}}),diameter(\widetilde{\gamma }_H^{
\overline{t}})\})+5.
\end{array}
$$
If we prove that $(\widetilde{f}_{\overline{t}})^N(\theta _{\overrightarrow{v^{\perp }}}^{
\overline{t}})$ cannot have a topologically transverse intersection with $%
\theta _{\overrightarrow{v^{\perp }}}^{\overline{t}}+(a,b),$ it clearly cannot
intersect $\theta _{\overrightarrow{v^{\perp }}}^{\overline{t}}+(c,d).$ So, if
there were such a topologically transverse intersection, as
$\theta _{\overrightarrow{v^{\perp }}}^{\overline{t}}+(a,b)$ is disjoint
from $\theta _{\overrightarrow{v^{\perp }}}^{\overline{t}}$ and arcs
of stable manifolds shrink under positive
iterates of $\widetilde{f}_{\overline{t}},$ there would be some $(a^{\prime
},b^{\prime }),(a^{\prime \prime },b^{\prime \prime })\in {\rm Z%
\negthinspace
\negthinspace Z^2,}$ such that $\widetilde{P}_{\overline{t}}+(a^{\prime
},b^{\prime })$ belongs to $\theta _{\overrightarrow{v^{\perp }}}^{\overline{%
t}}$ and its unstable manifold has a transverse intersection with the stable
manifold of $\widetilde{P}_{\overline{t}}+(a^{\prime \prime },b^{\prime
\prime })$ which belongs to $\theta _{\overrightarrow{v^{\perp }}%
}^{t^{\prime }}+(a,b).$ As $\theta _{\overrightarrow{v^{\perp }}}^{\overline{%
t}}$ is bounded in the direction of $\overrightarrow{v}$ by $3+2.\max
\{diameter(\widetilde{\gamma }_V^{\overline{t}}),diameter(\widetilde{\gamma }%
_H^{\overline{t}})\},$ we get that $(a^{\prime \prime }-a^{\prime
},b^{\prime \prime }-b^{\prime }).\overrightarrow{v}>0,$ so lemma 0
implies the existence of a rotation vector outside
$\rho (\widetilde{f}_{\overline{t}}),$ a contradiction.

Thus, from the continuity of $t\rightarrow \theta _{\overrightarrow{v^{\perp }}%
}^t$ and $t\rightarrow \widetilde{f}_t,$ there exists $t^{\prime }\in ]
\overline{t},t^{*}[$ such that $(\widetilde{f}_{t^{\prime }})^N(\theta _{
\overrightarrow{v^{\perp }}}^{t^{\prime }})$ has a non-topologically
transverse intersection with $\theta _{\overrightarrow{v^{\perp }}%
}^{t^{\prime }}+(c,d)$ and for all $t\in ]t^{\prime },t^{*}],$ the
intersection is topologically transverse. As above, from the fact that
stable manifolds shrink under positive iterates of $\widetilde{f}_{t^{\prime
}},$ the intersection that happens for $t=t^{\prime }$ corresponds to a
tangency between the unstable manifold of some translate of $\widetilde{P}%
_{t^{\prime }}$ which belongs to $\theta _{\overrightarrow{v^{\perp }}%
}^{t^{\prime }}$ with the stable manifold of some translate of $\widetilde{P}%
_{t^{\prime }}$ which belongs to $\theta _{\overrightarrow{v^{\perp }}%
}^{t^{\prime }}+(c,d).$ In other words, there is a tangency between $W^u(
\widetilde{P}_{t^{\prime }})$ and $W^s(\widetilde{P}_{t^{\prime
}})+(c^{*},d^{*})$ for some integer vector $(c^{*},d^{*})$ such that $\left|
(c^{*}-c,d^{*}-d).\overrightarrow{v}\right| \leq 3+2.\max \{diameter(
\widetilde{\gamma }_V^{\overline{t}}),diameter(\widetilde{\gamma }_H^{
\overline{t}})\}.$ This estimate follows from the fact that if $\widetilde{P}%
_t+(e,f)$ belongs to $\theta _{\overrightarrow{v^{\perp }}}^t,$ then $\left|
(e,f).\overrightarrow{v}\right| \leq 3+2.\max \{diameter(\widetilde{\gamma }%
_V^{\overline{t}}),diameter(\widetilde{\gamma }_H^{\overline{t}})\}.$ Here we are using that 
$\widetilde{P}_t \in \theta _{\overrightarrow{v^{\perp }}}^t.$ 

As $t^{*}>\overline{t}$ is arbitrary, if we remember that for a generic
family as we are considering, topologically transverse intersections are $%
C^1 $-transverse and tangencies are always quadratic, the proof of the main
theorem is almost complete. We are left to deal with the last part of the
statement, which says that for all parameters $t>\overline{t},$%
$$
W^u(\widetilde{P}_t)\text{ has transverse intersections with }W^s(\widetilde{%
P}_t)+(a,b),\text{ }\forall (a,b)\in {\rm Z\negthinspace
\negthinspace Z^2}. 
$$

As for $t>\overline{t},$ $(0,0)\in int(\rho (\widetilde{f}_t)),$ the main
result of \cite{c1epsilon} implies that (for each $t>\overline{t})$ $
\widetilde{f}_t$ has a hyperbolic periodic saddle (not necessarily fixed) $
\widetilde{Z}_t\in {\rm I\negthinspace R^2}$ such that $W^u(\widetilde{Z}_t)$
has transverse intersections with $W^s(\widetilde{Z}_t)+(a,b),$ $\forall
(a,b)\in {\rm Z\negthinspace
\negthinspace Z^2.}$ So, as $W^u(\widetilde{P}_t)$ and $W^s(\widetilde{P}_t)$
are both unbounded, exactly as we did in the proof of the fact from step 1
of this proof, $W^u(\widetilde{P}_t)$ has transverse intersections with $W^s(
\widetilde{Z}_t)$ and $W^s(\widetilde{P}_t)$ has transverse intersections
with $W^u(\widetilde{Z}_t).$ Thus an application of the $\lambda $-lemma
concludes the proof.

\vskip 0.2truecm

{\it Acknowledgments: }

When finishing the writing of this paper, I noticed that if I did not assume 
$C^1$-transversal intersections between stable and unstable manifolds of
hyperbolic periodic saddles, but only topologically transversal
intersections in certain parts of the arguments, the proofs would become
much more complicated. So, I started looking for a result saying that for
generic 1-parameter families of surface diffeomorphisms, homoclinic and
heteroclinic points are either quadratic tangencies or $C^1$-transverse. At
first, as I could not find it, I got quite anxious and asked several people
who have worked on the subject about this. Here I thank them for their
answers, which were positive, as this kind of result was known to be true
since \cite{nepata}. My thanks go to Carlos Gustavo Moreira, Eduardo Colli,
Enrique Pujals, Marcelo Viana and Pedro Duarte.

I also would like to thank Andres Koropecki for conversations on the results
of \cite{oliveira} and for explanations on his results in \cite{patmeykoro1}
and \cite{patmeykoro2}.

%
%
%
%

\end{document}